# Comparing dragonfly wings to jars of marbles through the lens of hyperuniformity

Karen E. Daniels, Department of Physics, North Carolina State University
Charles Emmett Maher, Department of Mathematics, University of North Carolina at Chapel Hill
Katherine Newhall, Department of Mathematics, University of North Carolina at Chapel Hill
Mason A. Porter, Department of Mathematics, University of California, Los Angeles; Department of Sociology, University of California, Los Angeles; Santa Fe Institute
Christopher Rock, Center for Additive Manufacturing and Logistics, North Carolina State University

## Abstract

When we look at the world around us, we see both organized (also called **ordered**) and disorganized (also called **disordered**) arrangements of things. Carefully-tiled floors and brick walls have organized and repeating patterns, but the stars in the sky and the trees in a forest look like they're arranged in a disordered way. We also see objects, like jars of marbles and the lacy wings of insects, that lie between ordered and disordered extremes. Although the marbles in a jar don't sit on a regular grid like carefully-arranged tiles, the collection of marbles does have some consistent features, such as the typical size and spacing between them. However, the positions of the marbles are much less random than the positions of the stars in the sky. To help understand and classify these patterns, mathematicians and physicists use the term **hyperuniform** to help them describe the situations of being perfectly organized or being disorganized in an organized way. In this article, we discuss various fascinating properties of hyperuniform patterns. We explore where they occur in the natural world and how engineers are using them to build new structures.

## What types of order do we see in the world around us?

Imagine that you want to give instructions to build a tiled floor from 10 cm × 10 cm square tiles. Once the first tile is placed, the next tile needs to be placed exactly 10 cm either to its left or right or to its top or bottom. And then all of the other tiles also need to be placed in this way until the floor is completely covered with tiles. It requires skill to do this well. However, when we dump a bunch of marbles into a jar, they all end up being about the same distance from each other without much effort. Amazingly, when we look at the world around us, we see many examples of this phenomenon. We show a few examples in Figure 1.

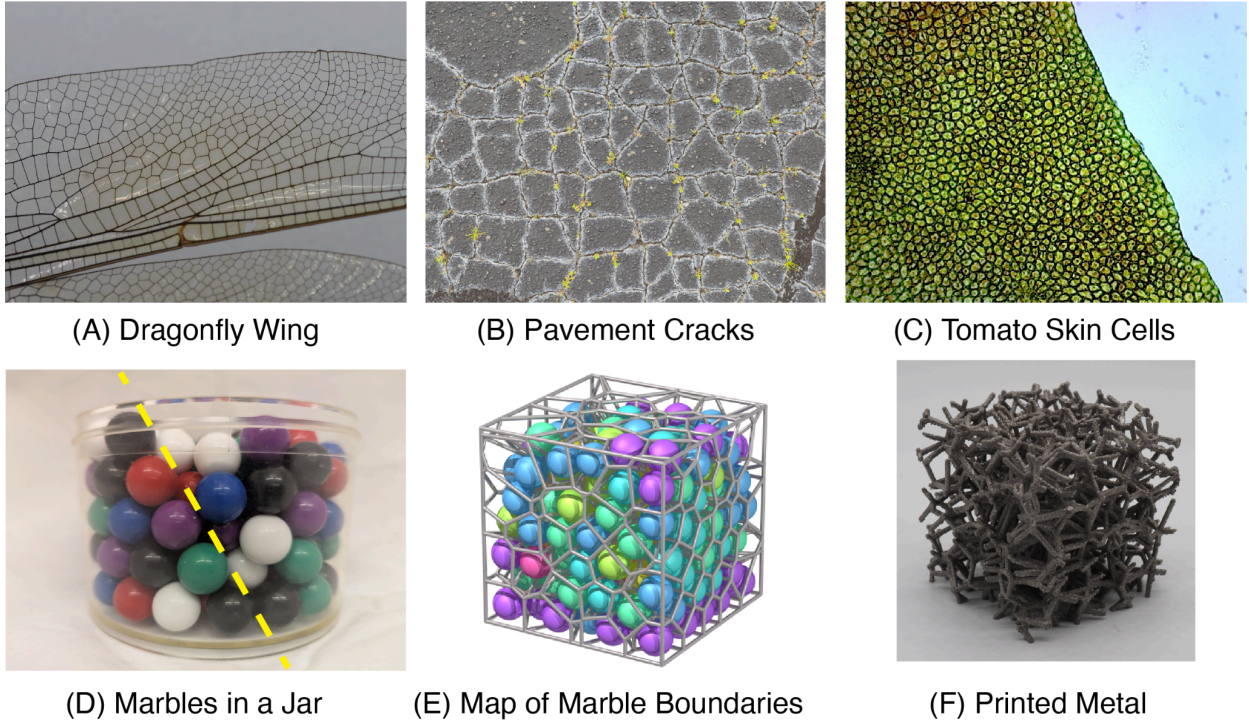

**Figure 1:** Six examples of disordered patterns. (A) A dragonfly wing from the collection of the NC State Insect Collection. More information about these patterns is available in [1]. (B) Cracks in the pavement of a parking lot. More information about these types of patterns is available in [2]. (C) Skin cells in a Roma tomato skin (Image credit: Wikimedia) (D) Marbles in a jar, with a more-ordered region on the left side of the dashed line and a less-ordered region on the right side of the line. For more information about packing spheres in a jar, see this video by Eugenia Cheng. (E) A digital blueprint showing how to draw boundaries between marbles (Image credit: Chris Rycroft, https://math.lbl.gov/voro++/examples/custom_output/). (F) A real-object version of a structure that is similar to the one in (E). [All photos are by the authors unless otherwise credited.]

One example is the lacy wings of dragonflies in Figure 1A. These wings have vein-like structures called "ribs", which form an imperfect grid that is similar to the patterns in cracked pavement (see Figure 1B)) and the walls of cells in the skin of a tomato (see Figure 1C). You can see that the cells have both a typical size and a typical number of neighboring cells, but you can *also* see that the depicted **tiling** is completely disorganized, unlike most floor tiles (which are often laid down carefully by hand). However, although the patterns in the tomato skin, the dragonfly wing, and the pavement are all disorganized, the wing has mostly empty space and the pavement has only thin regions of empty space within its cracks. By contrast, the tomato cells fill all of the space, with its cell walls composed of a different material than the stuff (which is called "cytoplasm") within its cells.

The examples in Figures 1A–C are all in a thin, flat layer, but we can also think about arrangements in three-dimensional (3D) space. When we put marbles into a clear box (see Figure 1D), we can either do it carefully (as in the organized pattern on the left side of the dashed line) or messily (as in the disorganized pattern on the right side of the dashed line). The two sides of the box share the property of having marbles separated by almost equal distances, but the direction (up, down, left, and right) between the marbles is consistent on the left side and inconsistent on the right. With a computer, we can write a program that draws the boundaries (see Figure 1E) between the space that "belongs" to each marble. The pattern on the sides of the cube looks a lot like the dragonfly wing. Once we've drawn these disordered bars, engineers can create them in real life (see Figure 1F), and then we can study their properties.

## How can we classify the patterns that we see?

Our eyes tell us that all of the objects in Figure 1 are disordered, but it would be nice to measure how disordered they are using numbers. Ideally, we'd even like to be able to determine which objects are more disordered than other objects.

To build ideas about how to do this, we can start by looking at the dots (also called **point patterns**) in Figure 2. The one in panel A is disordered and consists of randomly scattered points, while the one in panel C is completely ordered. But what about the point pattern in the middle? The spacing between the points is fairly even (somewhat like the pattern on the right), but the points have some randomness to their positions (like the pattern on the left).

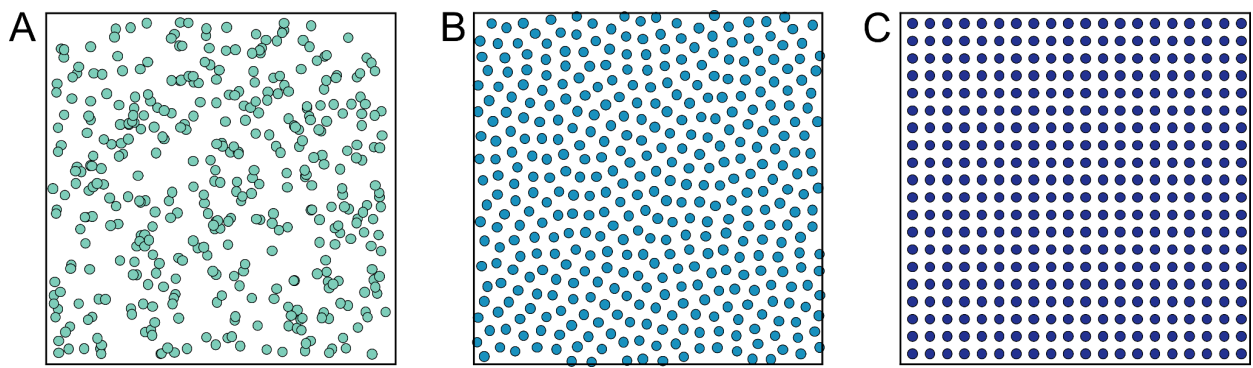

**Figure 2:** Point patterns with three different types of order/disorder. Panel A has a totally disordered point pattern, while panel C has a completely ordered point pattern. Panel B is between these extremes; this point pattern is a disordered hyperuniform point pattern with small-scale disorder (like panel A) and large-scale order (like panel C).

In Video 1, we illustrate one approach that we can use to answer this question. Let's start by drawing a circle inside the pattern and then counting the number of points that are inside it. We can then do this again and again for lots of different circles of the same size. Sometimes there may be 100 points, sometimes there may be 107 points, and sometimes perhaps there may be only 88 points. We measure how inconsistent these counts are, for different circle sizes, using a number called the **variance**. Researchers have discovered that measuring how fast this variance grows as we make the circle progressively larger (in other words, as we increase the **length scale** of the circle's radius**)** allows us to measure how disordered a point pattern is [3,4]. In this way, we can classify the presence and amount of disorder in a pattern by looking at the *trend* in the inconsistency of the counts for progressively larger length scales, rather than by choosing a single length scale when we compare values.

Now let's watch Video 1 again and pay attention to the inconsistency of the counts as the circle grows. For the most disordered point pattern (see Figure 2A), the variance grows more with circle size than it does for the perfectly ordered point pattern (see Figure 2C). The disordered point pattern in Figure 2B is like the pattern in Figure 2A in being disordered, but it's like the pattern in Figure 2C in how its variance grows. This special property — looking disordered, but doing so in an ordered way if you look through a large enough window, such as through a large enough circle — is what researchers call disordered **hyperuniformity**. Notably, we observe this special property because we track the trend of the variance as we change the circle sizes.

Scientists and engineers are very interested in objects that simultaneously share some properties of the disordered pattern in Figure 2A and some properties of the perfectly ordered pattern in Figure 2C. Testing for hyperuniformity with growing circles also works for any of the pictures in Figures 1A–C by calculating how much of each circle is dark and how much of each circle is light. Similarly, for the 3D objects in Figures 1D–F, we can grow a spherical bubble and calculate how much of each bubble is light and how much of it is dark.

## What different types of hyperuniformity are there?

Imagine that you are given a box of identical square tiles to retile a floor without leaving any large gaps. How many ways are there to do this? One popular method is to first place a straight line of tiles and then make sure that all subsequent lines match this pattern. Perhaps you also decide to offset the next line, as in a brick wall. However, what if you have circular tiles instead of square ones? You can still use both of these methods, but there are also *lots* of other options. If you have various coins or cookies of various sizes around, try this yourself. Place them randomly on a table and then use your hands to fit them into a smaller region. Every time you do this, you obtain a similar pattern as your previous ones, but the specific locations are different each time. The same is true if you put marbles into a jar. How many ways are there to do this? (You can read about tiles with even more shapes in this *FYM* article about different materials that things are made of.)

Should we count each version separately? Or should we instead count the whole group as one item because they all look similar and are generated using the same process? Researchers may choose to look at the exact properties of one version. (For example, how well does it conduct electricity?) However, they may study either the typical value (for example, the typical polygon area in Fig1A–C) or how many times the different-sided polygons appear in an entire collection of different versions that are made using the same recipe.

We can also change the recipe and create even more patterns, each of which is part of their own collection. Compare the different panels of Figure 3. Each panel shows a particular arrangement that is made using some set of rules, with a different set of rules for each panel. For each panel, do you think that the collection has one possible arrangement or many? What different shapes occur? Are the areas inside these shapes always similar, or can they be very different from each other?

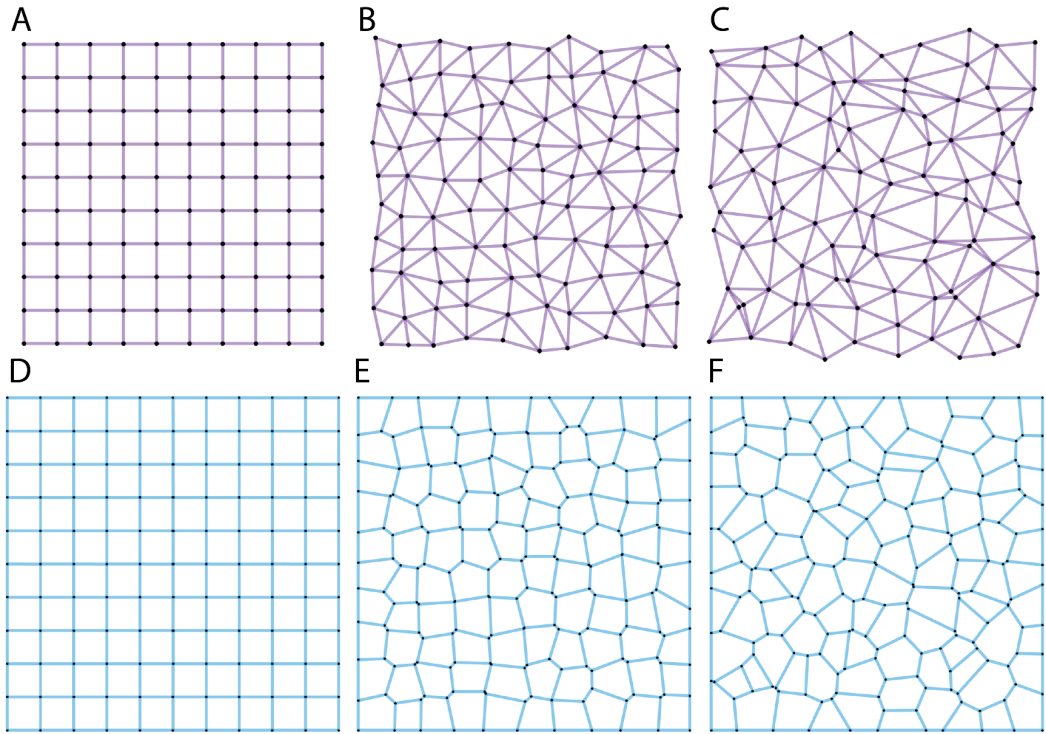

**Figure 3:** There are limited ways to be ordered, like the grids in panels A and D. Panels B and E show two different ways to be somewhat disordered, and panels C and F show two different ways to be even more disordered.

We made each of the patterns in Figure 3 starting from a point pattern that is similar to one in Figure 2. The top row uses a connection-based rule. We connect the random points to their nearest points in a way that ensures that no two lines overlap or cross. (The official name for this set of rules is a *Delaunay triangulation* [5].) The bottom row uses a polygon-based rule. The center of each polygon has a point from the point pattern, and the polygons equally divide the space that surrounds each point. (The official name for this set of rules is a *Voronoi tessellation* [5].) The first method (top row) tends to produce triangles, and the second method (bottom row) tends to produce polygons with both more sides and different numbers of sides.

The thing that distinguishes patterns in the different columns is that the initial point pattern is less ordered as we proceed from left to right. When we start from a perfect square grid of points (see Figures 3A,D), the two sets of rules agree on one perfect square tiling. However, when we start from a disordered point pattern, there are many more possible tilings. This produces a huge range of possibilities for engineers. With a new random point pattern, we obtain a new pattern of triangles or polygons (depending on which rule we follow). When we start from a very disordered point pattern (see Figures 3C,F), we obtain triangles or polygons that are all very different from each other. However, when we start from a hyperuniform point pattern, we obtain triangles or polygons whose sizes and shapes are much more similar to each other.

How can engineers take advantage of the large design space that disordered arrangements provide? One potentially interesting property is that materials that we produce with these arrangements don't have any special directions (like the up–down and right–left directions on a square grid). Therefore, if we try to run electricity or some other signal from one side of a disordered arrangement to the other, there is no direct route to take. Researchers are currently seeking ways to take advantage of this property.

# Conclusions and Outlook

We started off by asking how to make an unlikely comparison between dragonfly wings and jars of marbles. In both cases, the objects look disordered, and now we know some ways to measure this disorder. Take a look at the top row of Figure 1 again, and observe that the pictures in it look a bit like the pictures in the middle column of Figure 3! It turns out that we can also draw the jar of marbles in a similar way if we use a computer to map the boundaries between the marbles (see Figure 1E, and [6]. Now we see that the boundaries between the marbles look a lot like the ribs of a dragonfly wing. If we want to do a more detailed analysis, we can use the variance techniques to compare their amount of disorder. These types of ideas were used very recently in [7] to study the dots on the skin of squids!

We have also seen that patterns in natural and engineered systems can have different types of order, ranging from perfectly ordered rectangular bricks to completely disordered stars in the sky. Many natural systems have disorder levels between these extremes, and designable disorder can help engineers be more creative and successful. For example, if you're putting cellular-phone or electrical towers up in your town, everybody wants to be close enough to a tower to get good coverage, but nobody wants a tower on top of their house. Therefore, it may be useful to arrange the towers using a disordered grid. Additionally, a lot of sports equipment depends on different parts of it being stiff or flexible (to help protect athletes from impact injuries), while still being lightweight. Engineers are just getting started at developing applications of disordered and hyperuniform point patterns, and maybe you will develop the next great design.

Designers and engineers can use technology like 3D printing to create objects such as the one in Figure 1F or in [this FYM article](#). 3D printing, which is sometimes called additive manufacturing, allows us to create a real object by depositing material using the instructions in a digital blueprint, typically by adding one thin layer at a time. (Your local school or library may even have a 3D printer in its [Makerspace](#)!) These digital blueprints can either be designed (on a computer) with different types of disorder or inspired by natural, disordered structures like dragonfly wings. Sometimes, 3D printing software automatically creates hollow spaces within a solid object to save on material, and instead fills that space with open patterns like Figure 3. Researchers and engineers are currently investigating how they can use the type and amount of disorder in objects to enhance their function or even create new functions!

## Acknowledgements

We thank our young readers (Nia Chiou, Taryn Chiou, Zoe Chiou, Isaac Chiu, Jacob Chiu, Stan, Yarden, and Yuval) and their teachers and relatives (Lyndie Chiou, Christina Chow, and Veronika Juylova) for helpful comments. We thank the NC State Insect Collection for loaning us dragonflies to photograph. All authors were supported by a collaborative NSF DMREF grant, with grant numbers CMMI-2323341 (KED, CDR), CMMI-2323342 (KAN, CEM), and CMMI-2323343 (MAP).

## Glossary

- **Ordered (versus disordered):** A pattern that overlaps with itself (or doesn't) when you slide it in at least one direction.
- **Length Scale**: A distance or size that you are using for your observations or measurements.t
- **Hyperuniform:** A pattern that, on large length scales, behaves like it is ordered even when it looks disordered.
- **Point Pattern:** A set of dots in a square, cube, or other region of space.
- **Tiling**: A set of one or more geometric shapes (called "tiles") that cover a region, such as the surface of a floor. A tiling is also known as a "tessellation".
- **Variance:** A statistical measure of how far a set of values are spread out.

# Conflict of Interest

The authors declare that the research was conducted in the absence of any commercial or financial relationships that could be construed as a potential conflict of interest.

# Author biographies

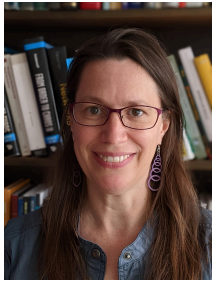

**Karen E. Daniels.** I am a physics professor now (at North Carolina State University), but decades of vacation photos reveal that I have always been interested in how patterns form in nature. After college, I got my start teaching middle and high school science, but I eventually went back for more school and received a PhD for my research doing laboratory experiments to study how patterns form in fluids. My lab and I often build experiments in the lab that teach us how the laws of physics operating on the small parts of a system can cause collective effects on the whole system.

**Charles Emmett Maher.** I am a postdoc in the Department of Mathematics at the University of North Carolina at Chapel Hill. I was born and raised near New Haven, Connecticut and have slowly been moving down the east coast of the United States since then. I study the structure of disordered systems using the tools talked about in this article, with a focus on dense collections of objects with many different shapes and sizes. I find these problems fun because despite how easy they are to say (for example, "How many marbles can I fit in this box?") we can learn a lot by thinking about different ways to solve them.

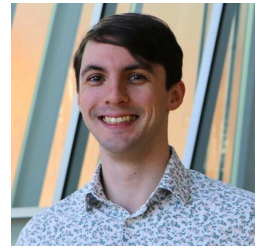

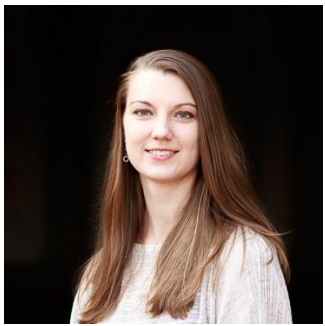

**Katherine (Katie) Newhall.** I am a professor at the University of North Carolina at Chapel Hill in the Department of Mathematics. I've been moving around the east coast, born in Maryland, schooling in New York (Troy and NYC), before settling in North Carolina. I study Stochastic processes, which are models that are disordered in time, for a diverse set of applications in physics and biology. Mathematics excites me because it's how I see the world around me: layers of simple rules for why objects or people move or behave as they do.

**Mason A. Porter.** I am a professor in the Department of Mathematics at UCLA. I was born in Los Angeles, California, and I am excited to be a professor in my hometown. I study networks, complex systems, nonlinear systems and their applications to physics, the social sciences, and many other areas. I am a big fan of games of all kinds, fantasy, baseball (Go Dodgers!), the 1980s, and other delightful things. I consider myself to be disorganized in an organized way.

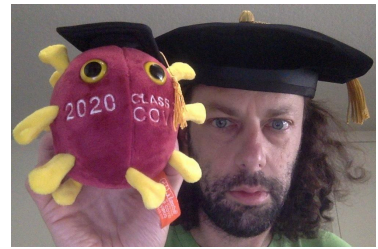

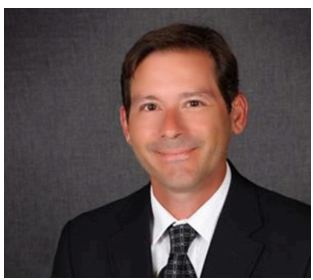

**Christopher Rock.** I am a research professor in the Center for Additive Manufacturing and Logistics at North Carolina State University. I am a life-long metallurgist and have manufactured metal shapes as large as tree trunks to smaller than a pinhead. Now I'm at NCSU CAMAL and use plastic and metal 3D printing to print 'math' and other really cool design ideas.